\documentclass[11pt]{amsart}   
\usepackage[greek,english]{babel}  
\usepackage{graphicx}
\usepackage{amsmath}
\usepackage{amssymb}
\usepackage{amsmath,amsthm}
\usepackage{verbatim}
\usepackage{sidecap}
\usepackage{tikz}
\usepackage{soul}
\newtheorem{theorem}{Theorem}[section]

\theoremstyle{definition}
\newtheorem{remark}[theorem]{Remark}

\newcommand{\R}{\mathbb{R}}

\newcommand{\Z}{\mathbb{Z}}

\titlepage 
\title[Parabolic solutions for the planar $N$-centre problem]{Parabolic solutions for the planar $N$-centre problem: multiplicity and scattering}
\author{Alberto Boscaggin, Walter Dambrosio and Duccio Papini}
\address{Alberto Boscaggin and Walter Dambrosio\newline \indent
 Dipartimento di Matematica ``Giuseppe Peano'', \newline \indent
Universit\`a di Torino, \newline \indent
Via Carlo Alberto, 10,
10123 Torino, Italy \newline \newline \indent
Duccio Papini \newline \indent
 Dipartimento di Scienze Matematiche, Informatiche e Fisiche, \newline \indent
Universit\`a di Udine, \newline \indent
Via delle Scienze, 206,
33100 Udine, Italy
}

\email{alberto.boscaggin@unito.it}
\email{walter.dambrosio@unito.it}
\email{duccio.papini@uniud.it}
\date{}

\begin{document}

\begin{abstract}
For the planar $N$-centre problem
$$
\ddot x = - \sum_{i=1}^N \frac{m_i (x-c_i)}{\vert x - c_i \vert^{\alpha+2}}, \qquad x \in \mathbb{R}^2 \setminus \{ c_1,\ldots,c_N \},
$$
where $m_i > 0$ for $i=1,\ldots,N$ and $\alpha \in [1,2)$, we prove the existence of entire parabolic trajectories, having prescribed asymptotic directions for $t \to \pm\infty$ and prescribed topological characterization with respect to the set of the centres.
\end{abstract}

%\date{\today}

\thanks{{\bf Acknowlegments.} Work partially supported by the 
ERC Advanced Grant 2013 n. 339958
{\it Complex Patterns for Strongly Interacting Dynamical Systems - COMPAT}, by the PRIN-2012-74FYK7 Grant {\it Variational and perturbative aspects of nonlinear differential problems} and by the INDAM-GNAMPA Project \textit{Dinamiche complesse per il problema degli $N$-centri}.}

\maketitle
\medbreak

\section{Introduction and statement of the main result}
\setcounter{page}{1}

The $N$-centre problem is the problem of the motion of a test particle in the attracting field generated by $N$ fixed heavy bodies $c_1,\ldots,c_N$; in Celestial Mechanics, it often arises as a simplified version of the restricted circular $(N+1)$-body problem in a rotating frame, when centrifugal and Coriolis' forces are neglected. For $N=1$, of course, it just reduces to the classical Kepler problem, while the case $N = 2$ has been solved by Jacobi (see, for instance, \cite{Whi59}). 
For $N \geq 3$, on the contrary, the problem has been proved to be analitically non-integrable \cite{Bol84} and,
in spite of its simple-looking structure, can indeed exhibit very complicated dynamics 
(see, among others, \cite{BolNeg03,BolNeg01,Dim10,KleKna92,Kna02,SoaTer12}). 
 
In this paper we will deal with the \emph{planar generalized $N$-centre problem}
\begin{equation}\label{eqmain}
\ddot x = - \sum_{i=1}^N \frac{m_i (x-c_i)}{\vert x - c_i \vert^{\alpha+2}}, \qquad x \in \mathbb{R}^2 \setminus \{ c_1,\ldots,c_N \},
\end{equation}
where $\alpha \in [1,2)$, thus including the classical Newtonian case $\alpha = 1$ as a particular case;
of course, $m_i > 0$ for $i=1,\ldots,N$.
Notice that the above equation has an Hamiltonian structure, with total energy given by
$$
H(x,\dot x) = \frac{1}{2} \vert \dot x \vert^2 - \sum_{i=1}^N \frac{m_i}{\alpha \vert x - c_i \vert^{\alpha}}.
$$
With this in mind, our aim is to prove the existence of entire \emph{parabolic} (i.e., zero-energy) solutions to \eqref{eqmain} having prescribed asymptotic directions at $\pm \infty$.
%these solutions are self-intersections free and separate the set 
%$\Sigma := \{c_1,\ldots,c_N\}$ of the centres according to any partition $\mathcal{P}$ of it. 
More precisely, denoting by $\Sigma = \{c_1,\ldots,c_N\}$ the set of the centres and naming partition of $\Sigma$ any subset $\mathcal{P} \subset \Sigma$ with $\mathcal{P} \neq \emptyset$ and $\mathcal{P} \neq \Sigma$, our main result reads as follows.

\begin{theorem}\label{thmain}
Let $N \geq 2$. For any asymptotic directions $\xi^-,\xi^+ \in \mathbb{S}^1$ with $\xi^- \neq \xi^+$ and for any partition 
$\mathcal{P}$ of $\Sigma$, there exists a self-intersection free parabolic solution $x: \mathbb{R} \to \mathbb{R}^2 \setminus \Sigma$ of \eqref{eqmain} satisfying $\vert x(t) \vert \to \infty$ for $t \to \pm \infty$,
$$
\lim_{t \to \pm \infty} \frac{x(t)}{\vert x(t) \vert} = \xi^\pm
$$
and separating the set $\Sigma$ according to the partition $\mathcal{P}$.
\end{theorem}
% (that is, $c_i \in \Gamma^1$ if  $c_i \in \mathcal{P}$ and $c_i \in \Gamma^2$ if $c_i \notin \mathcal{P}$, being $\Gamma^1,\Gamma^2$ the two connected components of $\mathbb{R}^2 \setminus x(\mathbb{R})$).
A comment about the statement: by the Jordan Theorem on a sphere, the above parabolic solution divides the plane into two connected components, both unbounded (see for instance \cite[Lemma 2.1]{BosGar16}); accordingly, the sentence 
``separating the set $\Sigma$ according to the partition $\mathcal{P}$'' means that two centres lie in the same connected component if and only if they are both in $\mathcal{P}$ or both in $\Sigma \setminus \mathcal{P}$.

Theorem \ref{thmain} has to be interpreted in the context of \emph{scattering}; indeed, it shows how the presence of two or more centres gives rise to (zero-energy) connections between any pair of asymptotic directions (but different), thus allowing in particular any value for the scattering angle. We stress that the analysis of the zero-energy case seems to be particularly interesting from this point of view; indeed, it is well known that for the central potential $V_\alpha(x) = \frac{m}{\alpha \vert x \vert^{\alpha}}$ (corresponding to the case $N=1$ in the generalized $N$-centre problem) all parabolic solutions span an angle of $2\pi/(2-\alpha)$ (see, for instance, \cite[Proposition 6.1]{BosDamTer17}). This is in strong contrast with the positive energy case, where all (but one) scattering angles are always achieved; accordingly, it is immediately understood that the possibility of an arbitrary zero-energy scattering angle as in Theorem \ref{thmain} is a genuine consequence of the presence of $N \geq 2$ centres and of the interaction of a parabolic solution with them. Incidentally, let us observe that, by collapsing all the centres into a single one, such parabolic solutions converge to the juxtaposition of two rectilinear zero-energy solutions of the $\alpha$-Kepler problem (see Remark \ref{eps}). From this perspective, we can also interpret Theorem \ref{thmain} as a continuation-type result, producing however classical solutions starting from generalized ones (the case $\xi^- = \xi^+$ being indeed the only one in which we cannot rule out the presence of collisions). 

We refer the reader to \cite{BarTerVer14,BarTerVer13,Che98,LuzMad14,MadVen09} for interesting investigations, from different point of views, about zero-energy solutions of various problems in Celestial Mechanics; we notice that, in spite of the differences between the considered models, all these results show the crucial role of parabolic solutions as carriers from different regions of the phase-space, in complete agreement with Theorem \ref{thmain}.  
We also mention that an extensive analysis of the scattering process for the planar $N$-centre problem has already be given in the excellent monograph \cite{KleKna92} by Klein and Knauf, dealing however only with the Newtonian case ($\alpha = 1$) and with positive energy solutions. The results therein are obtained via a global regularization of the problem, allowing to apply the theory of geodesics on surfaces of negative curvature. 
It is plausible that some results for the zero-energy case could be derived via a limiting procedure; we stress, however, that our approach is more direct and it allows the study of the generalized problem
\eqref{eqmain} with $\alpha \in [1,2)$ in a unified way. 

For the proof of our result, we combine indeed the variational approach to the construction of topologically non-trivial solutions of the Bolza boundary value problem associated with \eqref{eqmain}, developed in \cite{SoaTer12,SoaTer13}, together with a limiting procedure introduced in the recent paper \cite{BosDamTer17}, dealing with parabolic solutions of the $N$-centre problem in the three-dimensional space. Both these tools are available when $\alpha \in [1,2)$; it has to be emphasized, however, that the Newtonian case is still more difficult, and indeed requires the use of some (local, Levi-Civita type) regularization techniques. We also notice that, while in the spatial case solutions of the (fixed-energy) Bolza problem were found via a min-max argument, thus producing entire solutions with (at least generically) nontrivial Morse index, here minimization of the Maupertuis functional in suitable homotopy classes is enough, thus leading to locally minimal solutions. 

As a final comment, we remark that the multiplicity pattern in Theorem \ref{thmain} is a consequence of the result proved in \cite{SoaTer12,SoaTer13}, providing solutions separating the set of the centres according to any given partition of it. It is likely that the use of more refined arguments, on the lines of \cite{Cas17}, could lead to solutions in different homotopy classes, allowing for self-intersections and
revolutions around the centres; in this way, one should obtain a much richer zero-energy dynamics, including scattering solutions, semi-bounded solutions as well as bounded orbits exhibiting symbolic dynamics. All this will be the object of a future investigation.

\subsection{Plan of the paper.} In Section \ref{sez2} we review the existence of topologically non-trivial parabolic solutions of the Bolza problem, while in Section \ref{sez3} we show how to obtain entire parabolic solutions via a limiting procedure. Actually, we are going to prove that the conclusion of Theorem \ref{thmain} holds true for a larger class of equations of the type
\begin{equation}\label{eq-gen}
\ddot x = \nabla U(x), \qquad x \in \mathbb{R}^2 \setminus \Sigma,
\end{equation}
under suitable assumptions on the potential $U \in \mathcal{C}^\infty(\mathbb{R}^2 \setminus \Sigma)$ which we are going to list here below. First of all, we require 
\begin{equation} \label{eq-Upos}
U(x)>0,\quad \mbox{ for every } x\in \R^2\setminus \Sigma.
\end{equation}
Second, dealing with the behavior of $U$ near the centres we assume that, for some $\alpha \in [1,2)$,
\begin{equation}\label{vsing}
U(x) = \frac{m_i}{\alpha \vert x - c_i \vert^{\alpha}} + U_i(x), \qquad i=1,\ldots,N,
\end{equation}
where $m_i > 0$ and $U_i$ is smooth on $\mathbb{R}^2 \setminus (\Sigma \setminus \{c_i\}))$.
Finally, as for the behavior of $U$ at infinity, we require that, with the same $\alpha$ as above and some $m > 0$,
\begin{equation}\label{vinf}
U(x) = \frac{m}{\alpha \vert x \vert^{\alpha}} + W(x),
\end{equation}
where, for some $\beta > \alpha/2 + 1$,
$$
W(x) = O\left( \frac{1}{\vert x \vert^{\beta}}\right) \quad \mbox{ and } \quad
\nabla W(x) = O\left( \frac{1}{\vert x \vert^{\beta+1}}\right), \quad
\mbox{ for } \vert x \vert \to +\infty.
$$
It is easy to verify that the potential 
\begin{equation}\label{Nc}
V(x) = \sum_{i=1}^N \frac{m_i}{\alpha \vert x - c_i \vert^{\alpha}},
\end{equation}
giving rise to the generalized $N$-centre problem \eqref{eqmain}, satisfies all the above conditions, with $m = \sum_{i=1}^N m_i$ and $\beta = \alpha + 1$.

\section{Parabolic solutions of the Bolza problem}\label{sez2}

\noindent
In this section we look for solutions of the (free-time) fixed-endpoints problem
\begin{equation}\label{bol}
\left\{
\begin{array}{l}
\vspace{0.1cm}
\displaystyle{\ddot x = \nabla U(x)}\\
x(\pm \omega) = q^{\pm},
\end{array}
\right.
\end{equation}
saisfying the zero-energy relation
\begin{equation} \label{eq-energiaU}
\dfrac{1}{2}|\dot{x}|^2=U(x);
\end{equation}
recall that solutions of \eqref{bol} satisfying \eqref{eq-energiaU} are called parabolic solutions of \eqref{bol}.
Motivated by the final application, and in order to make all the discussion more transparent, we assume from the beginning that
\begin{equation}\label{hp-q}
\vert q^- \vert = \vert q^+ \vert \quad \mbox{ and } \quad q^- \neq q^+;
\end{equation}
also, we suppose that $\vert c_i \vert < \vert q^- \vert$ for $i=1,\ldots,N$, that is, all the centres lie inside the ball centered at the origin and of radius $\vert q^- \vert = \vert q^+ \vert$.

Having in mind a variational approach, we introduce the Maupertuis functional
$$
\mathcal{M}(u) = \int_{-1}^1 \vert \dot u(t) \vert^2 \,dt \int_{-1}^1 U(u(t))\,dt
$$
defined on the Hilbert manifold
$$
\widehat{\Gamma} = \widehat{\Gamma}_{q^{\pm}} = \Big\{ u \in H^1([-1,1];\mathbb{R}^2 \setminus \Sigma) \, : \, u( \pm 1) = q^{\pm} \Big\};
$$
notice that, in view of \eqref{eq-Upos}, it holds that $\mathcal{M}(u) \geq 0$ for any $u \in \widehat{\Gamma}$.
As well known (see, for instance, \cite[Theorem 4.1]{AmbCot93} and \cite[Appendix B]{SoaTer13}) $\mathcal{M}$ is smooth 
and any critical point $u \in \widehat{\Gamma}$ 
satisfies, for $t \in [-1,1]$, 
\begin{equation}\label{eqmau}
\ddot u(t) = \omega^2\, \nabla U(u(t)), \qquad \frac{1}{2}\vert \dot u(t) \vert^2 - \omega^2\, U(u(t)) = 0,
\end{equation}
where 
\begin{equation}\label{defomega}
\omega = \left(\frac{\int_{-1}^1 \vert \dot u(t)\vert^2 \,dt}{2 \int_{-1}^1 U(u(t))\,dt}\right)^{1/2}.
\end{equation}
Observe that, since $q^+ \neq q^-$, $u$ is not constant: as a consequence, $\omega> 0$ and the function
\begin{equation}\label{defxbeta}
x(t) = u\left( \frac{t}{\omega}\right), \qquad t \in [-\omega,\omega],
\end{equation}
is easily seen to be a parabolic solution of $\ddot x = \nabla U(x)$ on the interval $[-\omega,\omega]$; moreover, of course, $x(\pm \omega) = q^{\pm}$.

Following \cite{SoaTer12, SoaTer13}, multiple critical points of $\mathcal{M}$ can be found by minimizing in suitable homotopy classes. Precisely, write $q^{\pm} = \vert q^{\pm} \vert e^{i\theta_{\pm}}$, for suitable 
$\theta^{\pm} \in [0,2\pi)$, and define, for any $u \in \widetilde\Gamma$, the path $v_u: [-1,2] \to \mathbb{R}^2 \setminus \Sigma$ as
$$
v_u(t) = \begin{cases}
u(t) & t \in [-1,1] \\
\vert q^- \vert e^{i(\theta^+ + (\theta^- - \theta^+ + 2\pi)(t-1))} & t \in [1,2] 
\end{cases} \qquad \mbox{if } \theta^- < \theta^+,
$$
and
$$
v_u(t) = \begin{cases}
u(t) & t \in [-1,1] \\
\vert q^- \vert e^{i(\theta^+ + (\theta^- - \theta^+)(t-1))} & t \in [1,2] 
\end{cases} \qquad \mbox{if } \theta^+ < \theta^-,
$$
namely, we artificially close the path $u$ with the arc on $\partial B_{\vert q^- \vert}$ connecting $q^+$ with $q^-$ in the counterclockwise sense. With this notation, and given $l\in \Z_2^N$, we introduce the set
\[
\widehat{\Gamma}_l =\left\{u\in \widehat{\Gamma}:\ \textnormal{Ind}(v_u,c_i) \equiv l_i \mod 2,\ \forall \ i=1,\ldots, N\right\},
\]
being (in complex notation)
$$
\textnormal{Ind}(v_u,c_i) = \frac{1}{2\pi i} \int_{v_u} \frac{dz}{z - c_i}
$$
the usual winding number of a closed planar path. We are now in position to prove the following result:
\begin{theorem}\label{thfix}
Let $q^-, q^+$ be as in \eqref{hp-q} and let $l \in \Z_2^N$ satisfying
\begin{equation} \label{eq-lnonban}
\exists k\neq m: \; l_k\neq l_m.
\end{equation}
Then, there exists a self-intersection free parabolic solution of \eqref{bol}, corresponding to a (collision-free) minimizer of $\mathcal{M}$ in the $H^1$-weak closure of $\widehat{\Gamma}_l$.
\end{theorem}
\smallbreak
\begin{proof}[Sketch of the proof.] The existence of a minimizer $u$ of $\mathcal{M}$ in the $H^1$-weak closure of $\widehat{\Gamma}_l$ follows from standard lower-semicontinuity/coercivity arguments; notice however that the coercivity of $\mathcal{M}$ is not straightforward, following from the assumption at infinity \eqref{vinf} (see \cite[Lemma 4.2]{BosDamTer17} for the details). The fact that $u$ is collision-free can be proved as in \cite[Theorem 4.12]{SoaTer12} or in \cite[Theorem 2.3]{SoaTer13}, using \eqref{vsing} in an essential way and taking into account that the assumption $q^- \neq q^+$ rules out the case of collision-ejection solutions. Finally, the fact that $u$ is self-intersection free follows as in \cite[Theorem 4.12]{SoaTer12} again (see, in particular, \cite[Proposition 4.24]{SoaTer12}).
\end{proof}

%Following the discussion in \cite[pp. 3263-3264]{SoaTer12}, it is now possible to rephrase the above result in terms of partitions of the set of the centres.

%\beginningin{theorem}\label{thfix2}
%Let $q^-, q^+$ be as in \eqref{hp-q} and let $\mathcal{P}$ be a partition of $\Sigma$. 
%Then, there exists a self-intersection free parabolic solution of \eqref{bol} separating the set $\Sigma$ according to the partition $\mathcal{P}$.
%\end{theorem}

\section{Entire parabolic solutions}\label{sez3}

In this section we prove Theorem \ref{thmain} via an approximation argument. More precisely, given $\xi^-,\xi^+ \in \mathbb{S}^1$ with $\xi^- \neq \xi^+$ and a partition $\mathcal{P}$ of $\Sigma$, we first define $l \in \mathbb{Z}_2^N$ by setting 
$l_i = 1$ if and only if $c_i \in \mathcal{P}$ and we apply Theorem \ref{thfix} with the choice $q^{\pm} = R \xi^{\pm}$ for $R > 0$ large enough (notice that in this way \eqref{hp-q} is surely satisfied) so as to find an associated parabolic solution $x_R: [-\omega_R,\omega_R] \to \mathbb{R}^2$; then, we are going to show that an entire parabolic solution $x: \mathbb{R} \to \mathbb{R}^2$ can be obtained by passing to the limit when $R \to +\infty$. 

In order to do this, the assumption at infinity \eqref{vinf} will play a crucial role. For further convenience, we fix from the beginning two constants $C_-, C_+ > 0$ and a constant $K > \sup_i \vert c_i \vert + 1$ such that
\begin{equation}\label{stimaW}
\vert W(x) \vert \leq \frac{C_+}{\vert x \vert^\beta} \quad \mbox{ and }
\quad \vert \nabla W(x) \vert \leq \frac{C_+}{\vert x \vert^{\beta+1}}, \quad \mbox{ for every } \vert x \vert \geq K,
\end{equation}
\begin{equation}\label{stimaW2}
2 \vert W(x) \vert + \vert \nabla W(x) \cdot x \vert \leq \frac{(2-\alpha) m}{2 \alpha} \frac{1}{\vert x \vert^{\alpha}}, \quad \mbox{ for every } \vert x \vert \geq K,
\end{equation}
\begin{equation}\label{stimaV}
\frac{C_-}{\vert x \vert^{\alpha}} \leq U(x) \leq \frac{C_+}{\vert x \vert^{\alpha}}, \quad \mbox{ for every } \vert x \vert \geq K,
\end{equation}
and
\begin{equation}\label{stimak}
\sqrt{\frac{m}{\alpha}} \frac{1}{\vert x \vert^{\alpha/2}} - \frac{C_+}{\vert x \vert^{\beta-\alpha/2}} \leq 
\sqrt{U(x)} \leq
\sqrt{\frac{m}{\alpha}} \frac{1}{\vert x \vert^{\alpha/2}} + \frac{C_+}{\vert x \vert^{\beta-\alpha/2}}, \quad \mbox{ for every } \vert x \vert \geq K.
\end{equation}
The estimates \eqref{stimaW}, \eqref{stimaW2} and \eqref{stimaV} are rather obvious, while \eqref{stimak} follows from \eqref{stimaW} using the elementary inequalities
$1-\vert s \vert \leq \sqrt{1+s} \leq 1+\tfrac{1}{2}s$ (valid for $s \geq -1$).

We are now in position to give the proof; as a useful notation, we set $r_R(t) = \vert x_R(t) \vert$ and, whenever $r_R(t) \neq 0$, $s_R(t) = \tfrac{x_R(t)}{r_R(t)}$. We split our arguments into several steps; first of all, we observe that due to the assumption \eqref{eq-lnonban} any solution $x_R$ enters the ball $B_K$, so that 
$$\limsup_{R \to +\infty}\min_t r_R(t) \leq K < +\infty.
$$

\subsection{The virial identity and some preliminary estimates}\label{subsez3.1}

Preliminary, we observe that, due to the fact that $x_R$ has zero-energy (see \eqref{eq-energiaU}), the following equality 
- often referred to as virial identity - holds true:
\begin{equation}\label{vir}
\frac{d^2}{dt^2} \left(\frac{1}{2} r_R(t)^2\right) = 2 U(x_R(t)) + \nabla U(x_R(t)) \cdot x_R(t).
\end{equation}
Using \eqref{vinf} and \eqref{stimaW2},
we see that the above expression is strictly positive for $\vert x_R(t) \vert \geq K$, precisely
\begin{equation}\label{lagjac}
\frac{d^2}{dt^2} \left(\frac{1}{2} r_R(t)^2\right) \geq \frac{(2-\alpha)m}{2\alpha r_R(t)^{\alpha}}.
\end{equation} 
Therefore, $t_0 \in (-\omega_R,\omega_R)$ can be a local maximum for $t \mapsto r_R(t)$ only if
$r_R(t_0) < K$. 

On one hand, this implies that $r_R(t) < R$ for every $t \in (-\omega_R,\omega_R)$. As a consequence, $x_R$ separates the set $\Sigma$ according the partition $\mathcal{P}$, in the sense specified in \cite[pp. 3263-3264]{SoaTer12} (that is to say,
when closing the path $x_R$ as described in Section \ref{sez2} so as to find a Jordan curve $\gamma_R$, two centres lie in the same connected component of $\mathbb{R}^2 \setminus \gamma_R$ if and only if they are both in $\mathcal{P}$ or both in $\Sigma \setminus \mathcal{P}$).

On the other hand, it follows that there are exactly two instants $t^{\pm}_R \in (-\omega_R,\omega_R)$, with $t^-_R < t^+_R$, such that $r_R(t^\pm_R) = K$ (implying $r_R(t) < K$ for $t \in (t^-_R,t^+_R)$ and $r_R(t)> K$ for $t \notin [t^-_R,t^+_R]$); moreover, $\dot r_R(t) \neq 0$ for $t \notin (t^-_R,t^+_R)$. Using the fact that $x_R$ has zero-energy together with \eqref{stimaV}, we also find \begin{align*}
\omega_R - t^+_R & = \int_{t^+_R}^{\omega_R} \frac{\dot r_R(t)}{\dot r_R(t)}\,dt \geq \frac{1}{\sqrt{2C_+}} \int_{t^+_R}^{\omega_R} \frac{\dot r_R(t)}{r_R(t)^{-\alpha/2}}\,dt \\
& = \frac{1}{\left(1+\alpha/2 \right)\sqrt{2C_+}}\left(R^{1+\alpha/2} -
K^{1+\alpha/2} \right),
\end{align*}
implying that $\omega_R - t^+_R \to +\infty$ for $R \to +\infty$. Analogously, $-\omega_R - t^-_R \to -\infty$.
\smallbreak
For the rest of the proof, it is convenient to suppose $t^-_R = - t^+_R$, that is, the time spent by $x_R$ 
inside the ball $B_K$ is a symmetric interval with respect to the origin. This is not restrictive, up to a ($R$-dependent) time shift of the solution $x_R$. With a slight abuse of notation, we will still denote by $x_R$ this time-translation, and by 
$[\omega^-_R,\omega^+_R]$ its interval of definition.

\subsection{Passing to the limit: a generalized solution}\label{subsez3.2}

In this step, we show how to pass to the limit when $R \to +\infty$, so as to find an entire \emph{generalized solution}, that is, a solution with a zero-measure (but possibly non-empty) set of collision istants, see \cite{BahRab89}. For the next arguments, we write 
$$
\mathcal{A}_{[a,b]}(x) = \int_a^b \left( \frac{1}{2} \vert \dot x(t) \vert^2 + U(x(t)) \right)\,dt
$$
for the action of an $H^1$-path $x: [a,b] \to \mathbb{R}^2 \setminus \Sigma$; notice that, whenever $x$ satisfies the zero-energy relation \eqref{eq-energiaU}, we have 
$$
\mathcal{A}_{[a,b]}(x) = \int_a^b \vert \dot x(t) \vert^2 \,dt = 
2 \int_a^b U(x(t)) \,dt = \sqrt{2} \int_a^b \vert \dot x(t) \vert \sqrt{U(x(t))}\,dt.
$$
Having introduced this notation, the crucial point will be to prove that 
\begin{equation}\label{levest}
\limsup_{R \to +\infty}\mathcal{A}_{[t^-_R,t^+_R]}(x_R) < +\infty,
\end{equation}
with $t^{\pm}_R$ defined by the previous step. From this, several facts can be deduced. Precisely, since
$$
\left(\inf_{\vert x \vert \leq K} U(x) \right) \left( t^+_R - t^-_R\right) \leq \int_{t^-_R}^{t^+_R} U(x_R(t))\,dt, 
$$
we get at first that $t^+_R - t^-_R$ is bounded, say $t^+_R - t^-_R \leq 2T$ for any $R$. From this, together with the fact that $\vert x_R(t) \vert \leq K$ for $t \in [t^-_R,t^+_R]$ and with \eqref{levest} again, we infer that 
$$
\Vert x_R \Vert^2_{H^1(t^-_R,t^+_R)} = \int_{t^-_R}^{t^+_R} \left( \vert x_R(t) \vert^2 + \vert \dot x_R(t) \vert^2 \right)\,dt
$$
is bounded as well. Using moreover the fact that $\vert \ddot x_R(t) \vert \leq \left(\sup_{\vert x \vert \geq K} U(x)\right)$ for $t \notin [t^-_R,t^+_R]$, together with the boundedness of $\vert x_R(t^{\pm}_R) \vert$ and of $\vert \dot 
x_R(t^{\pm}_R) \vert = \sqrt{2 U(x_R(t^{\pm}_R))}$, we finally conclude that $x_R$ is bounded in $H^1_{\textnormal{loc}}(\mathbb{R})$.
As a consequence, there exists an $H^1$-function $x_\infty: \mathbb{R} \to \mathbb{R}^2$ such that $ x_R \to x_\infty$ weakly in $H^1_{\textnormal{loc}}(\mathbb{R})$ (in particular, uniformly on compact sets) for $R \to +\infty$. Of course, $x_\infty$ turns out to be a parabolic solution of \eqref{eq-gen} as long as it does
not collide with the set of the centres; moreover, $\vert x_\infty(t) \vert \geq K$ for $\vert t \vert \geq T$ so that the arguments of Subsection \ref{subsez3.1} imply that $x_\infty$ is unbounded for $t \to \pm \infty$.
Finally, by the $H^1_{\textnormal{loc}}$-boundedness and Fatou's lemma,
$$
\int_{-T}^{T} U(x_\infty(t))\,dt \leq \liminf_{R \to +\infty} \int_{-T}^{T} U(x_R(t))\,dt
= \liminf_{R \to +\infty} \frac12 \int_{-T}^{T} \vert \dot x_R(t) \vert^2 \,dt < \infty,
$$
implying that the set of collision instants has zero measure.

The rest of this subsection is then devoted to the proof of \eqref{levest}. We are going to show that
\begin{equation}\label{le1}
\mathcal{A}_{[\omega^-_R,\omega^-_R]}(x_R) \leq \left(\sqrt{\frac{2m}{\alpha}} \frac{4}{2-\alpha}\right) R^{1-\alpha/2} + M 
\end{equation}
and that
\begin{equation}\label{le2}
\mathcal{A}_{[\omega^-_R,t^-_R] \cup [t^+_R,\omega^+_R]}(x_R) \geq \left(\sqrt{\frac{2m}{\alpha}} \frac{4}{2-\alpha}\right) R^{1-\alpha/2} - M, 
\end{equation}
for some constant $M > 0$, from which \eqref{levest} clearly follows. 

We first prove \eqref{le1}. To this end, let us define the $H^1$-path
$$
\zeta(t) = 
\begin{cases}
\xi^+ \eta^+(t) & \; \mbox{ for } t \in [1,\Theta_R^+] \\
\gamma(s)(t) & \; \mbox{ for } t \in [-1,1] \\
\xi^- \eta^-(t) & \; \mbox{ for } t \in [\Theta^-_R,-1],
\end{cases}
$$
where $\gamma$ is an arbitrary $H^1$-path joining the points $K\xi^-$ and $K\xi^+$ and separating the set according to the partition $\mathcal{P}$ (in the sense specified in Section \ref{sez2}), $\eta^+:[1,+\infty) \to [K,+\infty)$ and $\eta^{-}: (-\infty,-1] \to [K,+\infty)$ are the solutions of the Cauchy problems
$$
\dot\eta^{\pm} = \pm\sqrt{2 U(\xi^{\pm}\eta^{\pm})}, \qquad \eta^{\pm}(\pm 1) = K
$$
and $\Theta^+_R,\Theta^-_R$ (for $R > K$) are the unique points such that $\eta^{\pm}(\Theta^{\pm}_R) = R$.
Then, we set
$$
\tilde \zeta(t) = \zeta\left( \Theta_R^- + \frac{1}{2} (\Theta^+_R - \Theta^-_R) (t+1)\right), \quad
\mbox{ for any } t \in [-1,1],
$$
in such a way that $\tilde \zeta$ is an $H^1$-path defined on $[-1,1]$, joining the points $R\xi^-$ and $R\xi^+$
and separating the set according to the partition $\mathcal{P}$. Using the well known relation
$$
\frac{1}{\sqrt{2}}\mathcal{A}_{[\omega^-_R,\omega^+_R]}(x_R) = \sqrt{\mathcal{M}(u_R)}, \qquad \mbox{ with } \, u_R(t) = 
x_R\left(\omega^-_R + \frac{1}{2}(\omega^+_R-\omega^-_R)(t+1)\right),
$$
together with the minimality of $u_R$ in the corresponding homotopy class, we find
$$
\frac{1}{\sqrt{2}} \mathcal{A}_{[\omega^-_R,\omega^+_R]}(x_R) \leq \sqrt{\mathcal{M}(\tilde \zeta)}.
$$
We therefore compute
\begin{align*}
\sqrt{\mathcal{M}(\tilde\zeta)} & = \frac{1}{\sqrt{2}} \inf_{\Theta > 0} \mathcal{A}_{[-\Theta,\Theta]}(\tilde \zeta(\cdot/\Theta))
\leq \frac{1}{\sqrt{2}}\int_{\Theta^-_R}^{\Theta^+_R} 
\left( \frac{1}{2} \vert \dot \zeta(t)\vert^2 + U(\zeta(t))\right)\,dt \\
& \leq 
\frac{1}{\sqrt{2}}\int_{\Theta^-_R}^{-1} 
\left( \frac{1}{2} \vert \dot \eta^-(t)\vert^2 + U(\xi^-\eta^-(t))\right)\,dt
+ \frac{1}{\sqrt{2}}\mathcal{A}_{[-1,1]}(\gamma)
\\
& \quad + 
\frac{1}{\sqrt{2}}\int_{1}^{\Theta^+_R} 
\left( \frac{1}{2} \vert \dot \eta^+(t)\vert^2 + U(\xi^+\eta^+(t))\right)\,dt \\
& = M_+ + \int_{-1}^{\Theta_R^-} \sqrt{U(\xi^-\eta^-(t))} \dot \eta^-(t) \,dt + \int_{1}^{\Theta_R^+} \sqrt{U(\xi^+\eta^+(t))} \dot \eta^+(t) \, dt
\\ & = M_+ + \int_{K}^R \sqrt{U(\xi^- r)}\,dr + \int_{K}^R \sqrt{U(\xi^+ r)}\,dr,
\end{align*}
with $M_+ = \tfrac{1}{\sqrt{2}}\mathcal{A}_{[-1,1]}(\gamma)$ (not depending on $R$).
Now, using the estimate from above in \eqref{stimak} we find
$$
\sqrt{U(\xi^{\pm}r)} \leq \sqrt{\frac{m}{\alpha}} \frac{1}{r^{\alpha/2}} + \frac{C_+}{r^{\beta - \alpha/2}}, \quad
\mbox{ for every } r \geq K,
$$
so that, with a simple computation,
$$
\sqrt{\mathcal{M}(\tilde\zeta)} \leq \left(\sqrt{\frac{m}{\alpha}}\frac{4}{2-\alpha}\right) R^{1-\alpha/2} + M_+ + \frac{4C_+}{2\beta - \alpha - 2},
$$
finally implying \eqref{le1}. To prove \eqref{le2}, we write
$$
\mathcal{A}_{[\omega^-_R,t^-_R] \cup [t^+_R,\omega^-_R]}(x_R) = \sqrt{2} \int_{[\omega^-_R,t^-_R] \cup [t^+_R,\omega^-_R]}
\vert \dot x_R(t) \vert \sqrt{U(x_R(t))}\,dt
$$
and we observe that $\vert \dot x_R(t) \vert \geq \vert \dot r_R(t) \vert$; moreover, by the arguments in Subsection \ref{subsez3.1}, $\dot r_R(t) < 0$ for $t \in [\omega^-_R,t^-_R]$ and $\dot r_R(t) > 0$ for $t \in [t^+_R,\omega^+_R]$.
Hence, using the estimate from below \eqref{stimak} yields the conclusion.

\subsection{Asymptotic directions}\label{subsez3.3}

We now prove that the (generalized) solution $x_\infty$ has $\xi^{\pm}$ has asymptotic directions for $t \to \pm\infty$, respectively; more precisely, writing $s_\infty(t) = \tfrac{x_\infty(t)}{\vert x_\infty(t) \vert}$ for 
$\vert t \vert \geq T$, we are going to show that $s_\infty(\pm \infty) = \xi^{\pm}$. Throughout this step of the proof, we assume that the solution $x_R$ is defined on the whole real line, as well. This is not restrictive, since the arguments of Subsection \ref{subsez3.1} (together with the boundedness of $\nabla U$ at infinity) rule out the occurrence of blow-up phenomena, and of course does not have influence on the local convergence $x_R \to x_\infty$; however, it turns out to be useful since it allows to perform estimates valid for any $t$ large enough (in absolute value).

We give the details for $t \to +\infty$. As a first step, we prove that
\begin{equation}\label{stima-r}
r_R(t) \geq \left( \frac{(2-\alpha)m}{2\alpha}\right)^{\frac{1}{\alpha+2}} \left( t - t^+_R\right)^{\frac{2}{\alpha+2}},
\quad \mbox{ for every } t \geq t^+_R.
\end{equation}
To obtain the above inequality, we first integrate \eqref{lagjac} on $[t^+_R,s]$, recalling that $ r_{R}(s) \le r_{R}(t) $
whenever $ t^{+}_{R} \le s \le t $, so as to obtain
$$
\frac{d}{dt} \left( \frac12 r_R(s)^2 \right) \geq \frac{(2-\alpha)m}{2\alpha m} \frac{(s-t^+_R)}{r_R(t)^\alpha}, \quad 
\mbox{ for every } t \geq s;
$$
a further integration on $[t^+_R,t]$ thus yields \eqref{stima-r}.

Taking into account that $t^+_R \leq T$, it follows from \eqref{stima-r} that there exists $\widehat{T} > T$ such that
$r_R(t) \geq K+1$ for $t \geq \widehat{T}$. We now claim that
\begin{equation}\label{stima-s}
\vert \dot s_R(t) \vert \leq \frac{C}{(t-T)^{\frac{4}{\alpha+2}}}, \quad \mbox{ for every } t \geq \widehat{T},
\end{equation}
where $C > 0$ is a suitable constant depending only on the potential (and on $K$).
To prove this, we define $A_R(t) = x_R(t) \wedge \dot x_R(t)$.
Taking into account \eqref{vinf} and \eqref{stimaW}, we first obtain from \eqref{stima-r} that
\[
| \dot A_R(t) |  = | x_{R}(t) \wedge \nabla W(x_{R}(t)) | \le \frac{ C_{+} }{ r_{R}^{\beta}(t) }
\le C_{+} \left[ \frac{ 2\alpha }{ (2-\alpha)m } \right]^{\beta/(\alpha+2)} \dfrac{1}{ (t-t^{+}_{R})^{2\beta/(\alpha+2)} }
\]
for every $ t \ge t^{+}_{R} $.
Denoting by $\hat t_R \in (t^+_R,\widehat{T}]$ the (unique) instant such that $r_R(\hat t_R) = K +1$, we then obtain,
for $t \geq \widehat{T}$,
\begin{align*}
\vert A_R(t) \vert & \leq \vert x_R(\hat t_R) \vert \vert \dot x_R(\hat t_R) \vert  +
\int_{\hat t_R}^{+\infty} \vert \dot A_R(\tau) \vert \,d\tau \\
& \leq r_R(\hat t_R) \sqrt{2U(x_R(\hat t_R))} + \int_{\hat t_R}^{+\infty} \vert \dot A_R(\tau) \vert \,d\tau,
\end{align*}
where
\begin{align*}
r_R(\hat t_R) \sqrt{2U(x_R(\hat t_R))} & \le \sqrt{2C_+} (K+1)^{\frac{2-\alpha}{2}} \\
\int_{\hat t_R}^{+\infty} \vert \dot A_R(\tau) \vert \,d\tau & \le C_{+} \left[ \frac{ 2\alpha }{ (2-\alpha)m } \right]^{\beta/(\alpha+2)}\frac{ \alpha+2 }{ 2\beta-\alpha-2 }
\dfrac{1}{ (\hat{t}_{R}-t^{+}_{R})^{(2\beta-\alpha-2)/(\alpha+2)} } \\
\hat{t}_{R} - t^{+}_{R} &\ge \dfrac{1}{ (1+\alpha/2)\sqrt{2C_{+}}} [ (K+1)^{1+\alpha/2} - K^{1+\alpha/2} ]
\end{align*}
using \eqref{stimaV} to bound from above $ U(x_R(\hat t_R)) $.
We have argued as in Subsection \ref{subsez3.1} to bound from below the quantity $\hat t_R - t^+_R$.
Observing that $\vert \dot s_R(t) \vert = \frac{\vert A_R(t) \vert}{\vert r_R(t) \vert^2}$ and using \eqref{stima-r} once again, \eqref{stima-s} finally follows.

From this we can easily conclude. Indeed, on one hand Lebesgue's theorem is
seen to apply, giving (together with uniform convergence on compact sets),
$$
s_R(+\infty) = s_R(\widehat{T}) + \int_{\widehat{T}}^\infty \dot s_R(\tau)\,d\tau \to s_\infty(\widehat{T}) + \int_{\widehat{T}}^\infty \dot s_\infty(\tau)\,d\tau = s_\infty(+\infty)
$$
for $R \to +\infty$. On the other hand, recalling that $s_R(\omega^+_R) = \xi^+$ and using \eqref{stima-s} again,
$$
s_R(+\infty) = \xi^+ + \int_{\omega^+_R}^{+\infty} \dot s_R(\tau)\,d\tau \to \xi^+, 
$$ 
finally yielding $s_\infty(+\infty) = \xi^+$. The proof that $s_\infty(-\infty) = \xi^-$ is analogous.

\subsection{Avoiding collisions}\label{subsez3.4}

In this step, we rule out the occurrence of collisions for $x_\infty$, that is, we prove that $x_\infty(t) \notin \Sigma$ for any $t \in \mathbb{R}$. We need to distinguish two cases, depending on whether $\alpha \in (1,2)$ or $\alpha = 1$.
\smallbreak
Let us suppose that $\alpha \in (1,2)$. Assume by contradiction that $x_\infty^{-1}(\Sigma) \neq \emptyset$; to fix the ideas, suppose that $x_\infty$ has (at least one) collision with the centre $c_1$ and take $\delta^* > 0$ so small that $c_i \notin B_{\delta^*}(c_1)$ for $i=2,\ldots,N$. Then
%, using the fact that the set of collision istants has zero measure (see Subsection \ref{subsez3.2} above), 
it is possible to find $\tau^-_R, \tau_R, 
\tau^+_R \in (t^-_R,t^+_R)$ such that 
$\tau^-_R < \tau_R < \tau^+_R$, $\delta_R := \vert x_R(\tau_R) - c_1 \vert = \min_t \vert x_R(t) - c_1 \vert \to 0^+$,
$$
\vert x_R(\tau^\pm_R) - c_1 \vert = \delta^*
\quad \mbox{ and }
\quad \vert x_R(t) - c_1 \vert \leq \delta^*, \quad \mbox{ for any } t \in [\tau^-_R,\tau^+_R].
$$   
Since $t^+_R - t^-_R$ is bounded and $x_R \to x_\infty$ uniformly on compact sets, both $\tau_R - \tau^-_R$ and $\tau^+_R - \tau_R$ are bounded away from zero. Let us define
$$
v_R(t) = \frac{1}{\delta_R}\left(x_R\left( \delta_R^{1+\alpha/2}t + \tau_R \right) - c_1 \right), \qquad t \in [-\gamma_R,\sigma_R],
$$
where
$$
-\gamma_R = \frac{\tau^-_R - \tau_R }{\delta_R ^{1+\alpha/2}} \quad 
\mbox{and} \quad  
\sigma_R = \frac{\tau^+_R - \tau_R }{\delta_R ^{1+\alpha/2}}.
$$
Notice that $-\gamma_R \to -\infty$ and $\sigma_R \to +\infty$, $\vert v_R(0) \vert = 1$, $\vert v_R(t) \vert \geq 1$ and $\vert \delta_R v_R(t)\vert \leq \delta^*$ for
$t \in  [-\gamma_R,\sigma_R]$. 
An easy computation shows that, writing $U$ as in \eqref{vsing}, $v_R$ satisfies
$$
\ddot v_R = - \frac{m_1 v_R}{\vert v_R \vert^{\alpha+2}} 
+\delta_R^{1+\alpha}\, \nabla U_1(\delta_R v_R + c_1)
$$
and
$$
\frac{1}{2}\vert \dot v_R \vert^2 = \frac{m_1}{\alpha \vert v_R \vert^{\alpha}} +  \delta_R^{\alpha} U_1(\delta_R v_R + c_1).
$$
As a consequence, it is easy to see
that $v_R \to v_\infty$ in $\mathcal{C}^2_{\textnormal{loc}}(\mathbb{R})$, where $v_\infty$ is a zero-energy solution of  
$$
\ddot v_\infty = - \frac{m_1 v_\infty}{\vert v_\infty \vert^{\alpha+2}}. 
$$ 
By \cite[Proposition 6.1]{BosDamTer17}, $v_\infty$ has transversal self-intersections. Since tranversal self-intersections are stable with respect to small perturbations, this contradicts the fact that $x_R$
(and hence $v_R$) is self-intersection free, thus ending the proof.
\smallbreak
Assume instead that $\alpha = 1$. Keeping the previous notation (and assuming now, up to passing to a subsequence, the existence of the limit $\tau_R \to \tau_\infty$) we define the Sundman integral
$$
s_R(t) = \int_{\tau_R}^t \frac{d\tau}{\vert x_R(\tau)-c_1 \vert}, \qquad t \in [\tau^-_R,\tau^+_R],
$$ 
and we use (with the usual identification $\mathbb{R}^2 \cong \mathbb{C}$) the well known Levi-Civita change of variables
$$
w_R(s)^2 = x_R(t_R(s)) - c_1, \qquad s \in [\sigma^-_R,\sigma^+_R], 
$$
being $t_R$ the inverse of $s_R$ and $\sigma^{\pm}_R = s_R(\tau^{\pm}_R)$.
Notice that the above change of variables is not one-to-one;
however, we can uniquely define $w_R$ by writing in polar coordinates $x_R - c_1 = \rho_R e^{i\varphi_R}$and setting
$w_R = \rho_R^{1/2} e^{i\varphi/2}$.
Also, observe that both $\sigma^-_R$ and $\sigma^+_R$ are bounded away from zero, since
$\vert s_R(t) \vert \geq \vert t - \tau_R \vert / \delta^*$. 

Standard computations yield:
\[
\ddot{x}_{R} = 2 w_{R}'' \dfrac{ w_{R} }{ |w_{R}|^{4} } = - 2 \frac{ w_{R}^{2}|w_{R}'|^{2} }{ |w_{R}|^{6} }.
\]
Here and in what follows all functions $ w_{R} $ and their derivatives are evaluated at $ s = s_{R}(t) $.
Using the equation and writing $U$ as in \eqref{vsing} with $\alpha = 1$ we get
\[
2 w_{R}'' w_{R} = 2 \frac{ w_{R}^{2}|w_{R}'|^{2} }{ |w_{R}|^{2} } - \frac{ m_{1} w_{R}^{2} }{ |w_{R}|^{2} }
+ |w_{R}|^{4}\nabla U_{1}( c_{1}+w_{R}^{2} )
\]
which gives
\[
w''_{R} = \frac{ w_{R} }{ |w_{R}| } |w_{R}'|^{2} - \frac{ m_{1} w_{R} }{ 2|w_{R}|^{2} }
+ \frac{ \overline{w}_{R} }{ 2 }|w_{R}|^{2}\nabla U_{1}( c_{1}+w_{R}^{2} )
\]
once it is multiplied by the complex conjugate $ \overline{w}_{R} $.
Finally the zero-energy relation for $x_R$ yields
\begin{equation}\label{lc-eq}
w_R'' = \frac{w_R}{2} U_1(c_1 + w_R^2) + \frac{\overline{w}_R}{2} \vert w_R \vert^2 \nabla U_1(c_1 + w_R^2);
\end{equation}
moreover
$$
\vert w_R(0) \vert = \vert x_R(\tau_R) - c_1 	\vert = \delta_R \to 0
$$
and
$$
\vert w_R'(0) \vert^2 = \frac{m_1}{2} + \frac{\vert w_R(0) \vert^2}{2}U_1(c_1 + w_R(0)^2) \to \frac{m_1}{2}.
$$
By a continuous dependence argument, $w_R$ converges (up to subsequences) uniformly on compact intervals containing the origin to the solution $w_\infty$ of the Cauchy problem associated with \eqref{lc-eq} having initial conditions $w_\infty(0) = 0$ and $w_\infty(0) = \nu$ for some $\vert \nu \vert^2 = \tfrac{m_1}{2}$; moreover, the symmetries of the differential equation \eqref{lc-eq} imply that it must be
$w_\infty(-s) = - w_\infty(s)$ for any $s$ small enough.

It follows that 
$$
t_R(s) = \tau_R + \int_0^s \vert w_R(\sigma) \vert^2 \,d\sigma \to t_\infty(s) = \tau_\infty + \int_0^s \vert w_\infty(\sigma) \vert^2 \,d\sigma
$$
uniformly on compact sets for $R \to +\infty$; moreover, the map $s \mapsto t_\infty(s) - \tau^*$ is an odd function. 
Taking into account that $x_R \to x_\infty$ uniformly on compact sets, we find
$$
w_\infty(s)^2 = x_\infty(t_\infty(s)) - c_1, \quad \mbox{ for every $s$ small enough}, 
$$
finally implying that $x_\infty(t) \neq c_1$ for $t$ near $\tau_\infty$ and that
$$
x_\infty(\tau_\infty - t) = x_\infty( \tau_\infty + t), \quad \mbox{ for every $t$ small enough.}
$$
Since $x_\infty$ is a classical solution of \eqref{eq-gen} outside the collision set, and possibly repeating the above argument for any collision instant, we find a contradiction with the global property that $x_\infty$ has different asymptotic directions for
$t \to \pm \infty$.

\subsection{Conclusion}\label{subsez3.5}

To conclude, we only need to show that $x_\infty$ is self-intersection free and that has the desired topological characterization. Actually, this second property immediately follows from the first one (taking into account the topological characterization of $x_R$), so let us show that $x_\infty$ is self-intersection free. Of course, transversal self-intersections are ruled out since $x_R$ is self-intersection free. On the other hand, assume by contradiction that there is a tangential self-intersection, that is, $x_\infty(t_1) = x_\infty(t_2)$ and $\dot x_\infty(t_1)$ parallel to $\dot x_\infty(t_2)$ for some $t_1 \neq t_2$. Then, the zero-energy condition gives $\vert \dot x_\infty(t_1) \vert = \vert \dot x_\infty(t_2) \vert$, so that
$\dot x_\infty(t_1) = \pm \dot x_\infty(t_2)$. Both the cases are not possible in view of the local uniqueness to the Cauchy problems: more precisely, in the first one $x_\infty$ should be periodic, while in the second one it should be
$x_\infty(t) = x_\infty(t_2+t_1-t)$, contradicting $\xi^- \neq \xi^+$.

\begin{remark}\label{asy}
Arguing as in \cite[Proposition 2.4]{BosDamTer17}, it is possible to prove that the above obtained parabolic solution
satisfies the asymptotic estimate
$$
\vert x_\infty(t) \vert \sim \left( \sqrt{\frac{m}{2\alpha}} (2+\alpha)\right)^{\frac{2}{2+\alpha}}\, \vert t \vert^{\frac{2}{2+\alpha}},
$$
when $t \to \pm\infty$.
\end{remark}

\begin{remark}\label{eps}
We finally briefly describe the behavior of the above found parabolic solutions when collapsing all the centres into a single one. In order to do this, we consider the parameter dependent $N$-centre problem
\begin{equation}\label{Nce}
\ddot y_\varepsilon = - \sum_{i=1}^N \frac{m_i (y-\varepsilon c_i)}{\vert y - \varepsilon c_i \vert^{\alpha+2}}
\end{equation}
when $\varepsilon \to 0^+$. Using a rescaling argument, solutions to the above equation can be obtained starting from solutions
of \eqref{eqmain}. More precisely, if $x_\infty$ denotes an entire parabolic solution of \eqref{eqmain}, then the function
\begin{equation}\label{defy}
y_\varepsilon(t) = \varepsilon \, x_\infty \left( \frac{t}{\varepsilon^{\frac{2+\alpha}{2}}}\right), \qquad t \in \mathbb{R};
\end{equation}
is a zero-energy solution of \eqref{Nce}. As a consequence of the asymptotic estimate given in Remark \ref{asy}, we have that the pointwise limit of $y_\varepsilon(t)$ for $\varepsilon \to 0^+$ exists, with
$$
\lim_{\varepsilon \to 0^+} y_\varepsilon(t) = \begin{cases}
\displaystyle \left( \sqrt{\frac{m}{2\alpha}} (2+\alpha)\right)^{\frac{2}{2+\alpha}} \vert t \vert^{\frac{2}{2+\alpha}}\xi^- & t < 0,
\\
0 & t = 0,
\\
\displaystyle \left( \sqrt{\frac{m}{2\alpha}} (2+\alpha)\right)^{\frac{2}{2+\alpha}} \vert t \vert^{\frac{2}{2+\alpha}}\xi^+ & t > 0.
\end{cases}
$$
As mentioned in the introduction, we have thus shown that, by collapsing all the centres into a single one, $y_\varepsilon$ converges to the juxtaposition of two rectilinear solutions of the $\alpha$-Kepler problem (actually, the convergence is easily seen to be $\mathcal{C}^2_{\textnormal{loc}}(\mathbb{R} \setminus \{0\})$; compare with \cite{Yu16}). 
\end{remark}

%\smallbreak
%\noindent
%{\bf Conflict of Interest.} The authors declare that they have no conflict of interest.

\end{document}